\documentclass[10pt,twosided]{article}
\usepackage{graphicx}
\usepackage{amsmath}
\usepackage{amssymb}
\usepackage{amscd}
\usepackage{hhline}
\usepackage{epic}
\usepackage{mathrsfs}

\newtheorem{theorem}{Theorem}

\newtheorem{theoremb}{Theorem}
\newtheorem{theoremc}{Theorem}

\newtheorem{theoreme}{Theorem}
\newtheorem{dfn}[theoremb]{Definition}
\newtheorem{rk}[theoremc]{Remark}

\newtheorem{lem}[theoreme]{Lemma}

\newtheorem{prop}[theorem]{Proposition}
\newenvironment{proof}[1][Proof]{\textbf{#1.} }{\qed \vspace{5pt}}
\newenvironment{Proof}[1]{\textbf{#1.} }

\newcommand\bib[1]{\bibitem[#1]{#1}}

\newcommand\abz{\hspace{13.5pt}}
\newcommand\qed{\phantom{\underline{y}}\hfill\hfill$\square$}
\newcommand{\comm}[1]{}

\newcommand\1{{\bf 1}}
\renewcommand\a{\alpha}
\renewcommand\b{\beta}

\newcommand\C{{\mathbb C}}

\renewcommand\d{\delta}

\newcommand\E{{\mathcal E}}

\renewcommand\l{\lambda}

\newcommand\La{\Lambda}

\newcommand\oo{\omega}
\newcommand\op[1]{\mathop{\rm #1}\nolimits}
\newcommand\ot{\otimes}
\newcommand\p{\partial}

\newcommand\R{{\mathbb R}}

\newcommand\ti{\tilde}

\newcommand\vp{\varphi}

\newcommand\x{\xi}

\makeatletter
\renewcommand{\@oddhead}{\hfil Differential Invariants \& PH-submanifolds\hfil}
\renewcommand{\@evenhead}{\hfil Boris Kruglikov \hfil}
\renewcommand{\@begintheorem}[2]{\begin{trivlist}\it
 \item[\hspace{\labelsep}{\bf #1\ #2.}]}
\renewcommand{\@endtheorem}{\end{trivlist}}
\makeatother

\newcommand\Cc{\let\mathcal\mathscr\mathcal C}

\begin{document}

 \title{Invariants and submanifolds\\ in almost complex geometry}
 \author{Boris Kruglikov}
 \date{}
% \date{\small Inst.\ of Math.\ and Stat., Univ.\ of Troms\o,
% Troms\o\ 90-37, Norway.\\
% E-mails: \quad kruglikov@math.uit.no, \quad lychagin@math.uit.no.
% \vspace{2pt}}
 \maketitle

 \vspace{-14.5pt}
 \begin{abstract}
In this paper we describe the algebra of differential invariants for $\op{GL}(n,\C)$-structures.
This leads to classification of almost complex structures of general positions.
The invariants are applied to the existence problem of higher-dimensional pseudoholomorphic
submanifolds.%
 \footnote{MSC numbers: 32Q60, 53C15; 53A55.\\ Keywords:
almost complex structure, equivalence, differential invariant, Nijenhuis tensor, pseudoholomorphic submanifold.}%
 \end{abstract}

%%%%%%%%%%%%%%%%%%%%%%%%%%%%%%%%%%%%%%%%%%%%%%%%%%%%%%%%%%%%%%%%%%%%%%%%%%%%
%0%
\section*{Introduction}

 \abz
Let $(M,J)$ be an almost complex manifold, $J^2=-\bf1$. In this paper we discuss only local aspects and so suppose $n=\frac12\dim M>1$. In this case the Nijenhuis tensor $N_J(\x,\eta)=[J\x,J\eta]-J[\x,J\eta]-J[J\x,\eta]-[\x,\eta]$ (which is a skew-symmetric (2,1)-tensor) is generically non-zero. Vanishing of $N_J$ is equivalent to local integrability of $J$ \cite{NN}.

It is known that all differential invariants can be expressed via the jet of the Nijenhuis tensor \cite{K$_1$}. In the first part of the paper we describe how this can be used to solve the equivalence problem of $\op{GL}(n,\C)$-structures. This problem is void for $n=1$ and was solved in \cite{K$_3$} for $n=2$, but we present here a uniform approach via differential invariants suitable for all $n$.

The differential invariants of an almost complex structure also occur in the problem of establishing pseudoholomorphic (PH) submanifolds. They played the crucial role in the proof of non-existence of PH-submanifolds for generic almost complex $(M,J)$ \cite{K$_2$}.

In dimension $2n\ge 8$ existence of a single higher-dimensional submanifold already imposes restrictions on the Nijenhuis tensor, so for their existence $N_J$ should be degenerate (though can stay far from being zero). On the other hand existence of 4-dimensional PH-submanifolds for 6-dimensional $(M,J)$ does not impose identities-restrictions on the tensor $N_J$ (there are open subsets of admissible tensors in the space of all Nijenhuis tensors).

At the second half of the paper we discuss when $(M,J)$ can have PH-foliations and when their number is bounded in non-integrable case.

%%%%%%%%%%%%%%%%%%%%%%%%%%%%%%%%%%%%%%%%%%%%%%%%%%%%%%%%%%%%%%%%%%%%%%%%%%%%
%1%
\section{First order invariants of almost complex structures}\label{S1}

 \abz
For $n=1$ the almost complex structures $J$ are complex and possess no local invariants.
So this case will not be considered in what follows.

For $n=2$ there are non-integrable structures, but there are no first order differential invariants. To explain this let us note that all such invariants must be derived from the Nijenhuis tensor $N_J$. In dimension 4 the linear Nijenhuis tensor (a purely tensorial object at a point, i.e. an element of $\La^2T^*\ot_{\bar\C}T$ with $T=T_xM$ \cite{K$_1$}) is special and the $\op{GL}(2,\C)\subset\op{GL}(4)$ orbit space consists of two points: zero and non-zero tensor $N_J$.

For non-zero tensor we can talk of the image $\Pi^2=\op{Im}(N_J)$ which, if we vary the point $x$, is a two-distribution in $M^4$, called Nijenhuis tensor characteristic distribution \cite{K$_2$}. Provided $J$ is generic, $\Pi^2$ is generic as well. In particular, there's the derived rank 3 distribution $\Pi^3=\p\Pi^2$. This leads to the fact that there's no second order invariants as well. However, we can associate the second order e-structure $\{\xi_i\}_{i=1}^4$, $\xi_i\in\mathcal{D}_M$, to $J$ as follows:
 $$
\x_1\in C^\infty(\Pi^2),\ \x_3\in C^\infty(\Pi^3),\ N_J(\x_1,\x_3)=\x_1,\ \x_2=J\x_1,\ [\x_1,\x_2]=\x_3,\ \x_4=J\x_3.
 $$
This defines $\xi_1,\x_2$ canonically up to $\pm1$ and $\x_3,\x_4$ absolutely canonically \cite{K$_3$}.

When $n=3$ there are moduli in the space of linear Nijenhuis tensors \cite{K$_2$}.
This is clearly seen from Theorem 7 loc.sit. Indeed the statement means that the space of differential invariants of order 1 is two-dimensional and the constants of the normal forms provide the invariants. In \cite{B} Bryant arrived independently to the result about dimension 2, observing that codimension of generic orbits w.r.t. the $\op{GL}(3,\C)\subset\op{GL}(6)$ action on the space of Nijenhuis tensors $\La^2T^*\ot_{\bar\C}T$ is 2 (where $T=T_xM$ is a model tangent space of dimension 6), because the stabilizer is two-dimensional.

Moreover in \cite{B} Bryant introduced some invariants of almost complex 6-dimensional  manifolds. All of them are expressed via a (1,1)-form $\oo$, which is given in coordinates via the components of the Nijenhuis tensor as follows:
 $$
\oo_i^j=\frac{N_{ik}^l\bar{N}_{jl}^k-N_{jk}^l\bar{N}_{il}^k}{\sqrt{-1}}.
 $$
Here complex coordinates adapted at the point to $J$ are used (in fact, in \cite{B} non-holonomic, i.e. frames). Note that this is a real-valued form and it can be written in invariant terms as follows (now we assume all tensors real):
 $$
\oo(\xi,\eta)=\op{Tr}[N_J(\xi,JN_J(\eta,\cdot))-N_J(\eta,JN_J(\x,\cdot))].
 $$
In particular, $\oo(\x,J\x)=2\op{Tr}[N_J(\x,N_J(\x,\cdot))]$ is not identically zero.

The form $\oo$ is $J$-compatible: $\oo(J\x,J\eta)=\oo(\x,\eta)$ and we can associate the quadric $q(\xi,\eta)=\oo(\x,J\eta)$, which equals
 $$
q(\xi,\eta)=\op{Tr}[N_J(\xi,N_J(\eta,\cdot))+N_J(\eta,N_J(\x,\cdot))].
 $$

Indeed, these both 2-tensors are skew-symmetric and symmetric parts of the form $T_{ij}=N_{ik}^l\bar{N}_{jl}^k$ and the pair $(S,\oo)$ forms a Hermitian metric provided $q$ (or equivalently $\oo$) is non-degenerate (it can be indefinite).

Let us investigate the form $\oo$ in terms of the normal forms given in \cite{K$_3$}:

 \begin{prop}
The (1,1)-form $\oo$ is degenerate precisely in the following cases in terms of the differential invariants from classification theorem 7 \cite{K$_1$}:
 \begin{itemize}
 \item[] \hspace{-15pt}{\rm NDG$_1$:} $\l=\pm1$, $\vp=0,\pi$.
 \item[] \hspace{-15pt}{\rm NDG$_2$:} $\vp=0,\pi$.
 \item[] \hspace{-15pt}{\rm NDG$_3$:} $\psi=\pm\frac\pi4\pm\frac\pi2=\pm\vp\pm\pi$.
 \item[] \hspace{-15pt}{\rm NDG$_4$,\ DG$_1$(4):}  Never.
 \item[] \hspace{-15pt}{\rm DG$_1$(1-3,5),\ DG$_2$(1-2):} Always.
 \end{itemize}
 \end{prop}

This is a straightforward tedious calculation. It shows generic non-degeneracy of the 2-form $\oo$. In \cite{B} Bryant discusses global implications of non-degeneracy, while the local aspects here show to which open strata should non-degenerate forms $\oo$ belong, which also yields topological restrictions on existence.

The canonical $G_2$-invariant almost complex structure  $J$ on $S^6$ corresponds to NDG.3 $\vp=0,\psi=\frac\pi2$, so in this case $\oo$ is non-degenerate. Also note that when the form $\oo$ is degenerate, then $M$ possesses a canonical distribution (kernel), which can be used to construct classification in the case of non-general position.

In dimension $n>3$ the orbit space of $\op{GL}(n,\C)$-action on $\La^2T^*\ot_{\bar\C}T$ is quite complicated. And indeed the space of invariants is pretty big, as will be discussed below.

%%%%%%%%%%%%%%%%%%%%%%%%%%%%%%%%%%%%%%%%%%%%%%%%%%%%%%%%%%%%%%%%%%%%%%%%%%%%
%2%
\section{General background on differential invariants}

 \abz
The equivalence problem of geometric structures on $M$ is usually solved either via differential invariants algebra or by constructing a canonical e-structure. In the first case the algebra can be represented either via some basic invariants and invariant differentiations or via some more differential invariants and Tresse derivatives.

However in the case of geometric structures the number of required differential invariants is smaller and equals $n=\dim M$, provided the restrictions of them to the structure are functionally independent (this is generically so).

Indeed, let $\pi$ be a bundle of geometric structures (associated with a tensorial bundle over $M$) and $\E$ a section of it (i.e. a geometric structure of the specified type), which can be represented as the image of a section $j:M\to E_\pi$. Let $\rho$ be the induced action of the pseudogroup $\op{Diff}_\text{loc}(M)$ on $\pi$ and $I$ be a differential invariant. Its restriction to $\E$ is the function $I_\E=j^*I\in C^\infty_\text{loc}(M)$.

Given $n$ functionally independent invariants $I^1,\dots,I^n$ we assume their restrictions $I^1_\E,\dots,I^n_\E$ are functionally independent (here and in what follows one can assume local treatment), so that they can be considered as local coordinates. Then one gets local frames $\p_i=\frac\p{\p I^i_\E}$ and coframes $\oo^i=dI^i_\E$. Any tensorial field $T$ can be expressed as $T=T_{i_1\dots i_s}^{j_1\dots j_t}\p_{j_1}\ot\dots\p_{j_t}\ot\oo^{i_1}\ot\dots\oo^{j_s}$ and the coefficients are scalar differential invariants.

Being expressed via $I^i_\E$ they form the complete set of invariant relations for equivalence problem. This is the principle of $n$-invariants \cite{ALV}.

Two remarks are of order. First: It is clear in this case that canonical frame field $\oo^i$ gives e-structure; otherwise around is also true, so that e-structures approach \cite{Ko} is equivalent to one with differential invariants. Second: Lifts of the derivations $\p_i$ are invariant differentiations and coefficients of $dJ|_\E=\frac{D J}{D I^i_\E}\oo^i$ are exactly Tresse derivatives of a differential invariant $J$ by the basis $I^i$ (\cite{Ku,KL}). Thus all the discussed approaches are equivalent.

Let us apply this to classification of almost complex structures of general position. This means that $\pi$ is the bundle of almost complex structures over $M$:
 $$
\pi^{-1}(x)=\{J\in\op{GL}(T_xM): J^2=-\1\}\simeq\op{GL}(2n,\R)/\op{GL}(n,\C)
\stackrel{\text{def}}=\mathcal{J}(2n).
 $$
The pseudogroup $G=\op{Diff}_\text{loc}(M)$ acts on $\pi$. The groupoid of its jets is denoted by $G^l\subset J^l(M,M)$, with natural projections being denoted by $\rho_{l,k}:G^l\to G^k$ and $\rho_l:G^l\to M$. Denote $G^l_x$ the fiber over the point $(x,x)\in G^0=M\times M$, which is also a sub-groupoid of $G$ called the differential group of order $l$ (we will sometimes omit reference to the point $x$).

Then $G^l_x$ acts on the fiber of $\pi_k:J^k\pi\to M$ over point $x$. Moreover denoting $\mathfrak{G}^{l+1}_x=\op{Ker}[\rho_{l+1,l}:G^{l+1}_x\to G^l_x]$ we obtain action of this normal subgroup on the fiber of the bundle $\pi_{l,l-1}$. Since for $l>0$ the group $\mathfrak{G}^{l+1}$ is abelian, the orbits in $F_l=\pi_{l,l-1}^{-1}(x_{l-1})$ are affine and so the differential invariants can be chosen affine in derivatives of order $l$. Usually they are non-linear in lower-order derivatives.

%%%%%%%%%%%%%%%%%%%%%%%%%%%%%%%%%%%%%%%%%%%%%%%%%%%%%%%%%%%%%%%%%%%%%%%%%%%%
%3%
\section{Equivalence problem for almost complex structures}

 \abz
Let $(M,J)$ be an almost complex manifold. If it is in general position, then as we have noticed above, it is enough to find $2n=\dim M$ differential invariants for local classification. Solution of the equivalence problem depends on $n$ (which we can assume to be $>1$).

 \begin{theorem}
The basic scalar differential invariants of $J$ solving the equivalence problem via the described methods can be specified as follows.
 \begin{itemize}
 \item[] \hspace{-15pt}$n=2:$ There are no differential invariants of order $\le 2$, but in order 3 there are (no less than) 4 differential invariants;
 \item[] \hspace{-15pt}$n=3:$ There are precisely 2 differential invariants of order 1 and 4 invariants of order 2;
 \item[] \hspace{-15pt}$n>3:$ There are at least $n^2(n-3)>2n$ differential invariants of order 1.
 \end{itemize}
 \end{theorem}

 \begin{proof}
Consider the cases.

${\bold n=2.}$ It is clear from the description in Section \ref{S1} that $J$ has no differential invariants of order 1 or 2. To get invariants of order 3 one proceeds as follows: the Maurer-Cartan coefficients $c_{ij}^k$ for the described canonical e-structure $\xi_s$ (defined by the 2-jet of $J$) are the invariants of order 3: $[\xi_i,\xi_j]=c_{ij}^k\xi_k$. Since $c_{12}^k=\d_3^k$ and $[\x_2,\x_4]$ can be expressed via other brackets from $N_J(\x_1,\x_3)=\x_3$, the number of such differential invariants is 16. Note that the invariants of \cite{K$_1$} \S6.1 (canonical 1- and 2-forms on $\Pi^2$) can be expressed via $c_{ij}^k$.

Let us notice that the result can be obtained via pure dimensional count. Indeed, rank of $\rho_{1,0}$ is 16 and that of $\pi$ is 8. The action of $\mathfrak{G}^1$ on $F_0=\pi^{-1}(x)$ is transitive (8-dimensional stabilizer). Rank of $\rho_{2,1}$ is 40 and that of $\pi_{1,0}$ is 32. Again the action of $\mathfrak{G}^2$ on $F_1$ is transitive (8-dimensional stabilizer). Next the rank of $\rho_{3,2}$ is 80 and that of $\pi_{2,1}$ is 80 as well, the corresponding action is transitive. Finally rank of $\rho_{4,3}$ is 140 and that of $\pi_{3,2}$ is 160. The action of $\mathfrak{G}^4$ cannot be transitive. Moreover if we consider the action of $G^4$ on $\pi_3^{-1}(x)$, it cannot be transitive as well, because even though the action of $G^3$ has $8+8=16$-dimensional stabilizer, the difference in dimension is $160-140-16=4$. Thus there are at least 4 differential invariants. Note though that there are more (as we explained above), so that the action of $G^4$ has a large stabilizer.

${\bold n=3.}$ We do the dimensional count. Rank of $\rho_{1,0}$ is 36 and that of $\pi$ is 18. The action of $\mathfrak{G}^1$ on $F_0$ is transitive (18-dimensional stabilizer). Rank of $\rho_{2,1}$ is 126 and that of $\pi_{1,0}$ is 108. It seems that the stabilizer should be 18-dimensional, but as we explained in Section \ref{S1} the action of
$G^2$ has orbits of codimension 2, so the dimension of stabilizer is by two bigger than can be expected. Next rank of $\rho_{3,2}$ is 336 and that of $\pi_{2,1}$ is 378, so that the pure difference of dimensions gives at least $378-336-18\cdot2-2=4$ differential invariants of order 3. All these invariants can be expressed via the normal forms of \cite{K$_3$}.

${\bold n>3.}$ Here the dimensional count can be misleading, so we better calculate codimension of orbits of $\op{GL}(2n)$-action on the space of linear Nijenhuis tensors. The stabilizer of a linear complex structure $J_0$ on $T$ is $\op{GL}(n,\C)$. Since $\dim\op{GL}(n,\C)=2n^2$ and $\dim\La^2T^*\ot_{\bar\C}T=n^2(n-1)$ the largest orbits have codimension $n^2(n-3)+\dim\op{St}$, where $\op{St}$ is a stabilizer of a generic point. The result follows.
 \end{proof}

This solves the equivalence problem for almost complex structures.

%%%%%%%%%%%%%%%%%%%%%%%%%%%%%%%%%%%%%%%%%%%%%%%%%%%%%%%%%%%%%%%%%%%%%%%%%%%%
%4%
\section{Existence of almost complex submanifolds}

 \abz
In a private communication M. Gromov asked the following question: how many higher-dimensional PH-submanifolds can an almost complex manifold possess? According to \cite{Gr,K$_2$} generically there are none.

On the other end, for integrable $J$ there're plenty. What happens in between? This question is quite difficult if PH-submanifolds are isolated, so we treat the case when they come in families, regulary fashioned, namely as PH-foliations.

We will consider in details 6-dimensional situation, the general case allows certain generalizations. Let us start with some examples.

 {\bf Example 1.}
Let $M=\C^3$ with almost complex structure $J$ being given in $2\times 2$ block form $J=\op{diag}(A_1,A_2,A_3)$, where the coefficients of $A_i\in C^\infty(M,\mathcal{J}(2))$ do depend on all 6 coordinates $(x^1,\dots,x^6)\in M$ in a generic way. Then the Nijenhuis tensor $N_J$ is non-degenerate. Indeed, we have for $i,j$ odd: $N_J(\p_i,\p_j)\in\C\langle\p_i\rangle\oplus\C\langle\p_j\rangle$, whence existence of kernel $\xi=\sum\limits_{i\text{ odd}}\a_i\p_i$ ($\a_i\in\C$ and multiplication means $\a\cdot\eta=\op{Re}\a\cdot\eta+\op{Im}\a\cdot{J}\eta$) implies that some of the vectors $N_J(\p_i,\p_j)$ have zero components in the above $\C^2$ decomposition.

In other words if we denote the above splitting as $M=V_1\oplus V_2\oplus V_3$, then the genericity condition is $N_J(V_i,V_j)\not\subset V_i$ for all $i\ne j$.

Thus we see that it is possible to have 3 transversal PH-foliation of $(M^6,J)$ with $J$ being maximally non-degenerate at each point. Generalization to dimension $2n$ is straightforward.

 {\bf Example 2.}
The above example can be modified as follows. Let $V_{ij}=V_i\oplus V_j$ and let the almost complex structure have a block form $J=\op{diag}(A,B)$ in the splitting $M=V_1\oplus V_{23}$, where $A\in C^\infty(M,\mathcal{J}(2))$ and $B\in C^\infty(M,\mathcal{J}(4))$. While $A$-block is allowed to be arbitrary, the $B$-block is assumed symmetric in $V_1$-direction, i.e. independent of $(x^1,x^2)$-coordinates.

Then any PH-curve $\mathcal{C}^2\subset V_{23}$ lifts to the 4D
PH-submanifold $V_1\times\mathcal{C}^2\subset M$. Thus we have an
infinite-dimensional family of PH-submanifolds
$\C\times\mathcal{C}^2$, all of which intersect by a leaf of the 2D
PH-foliation $V_1$.

Note that this family will persist if we allow the structure to have
the form
 $$
J=\begin{pmatrix} A & B \\ 0 & D\end{pmatrix}
 $$
in the splitting $M=V_1\oplus V_{23}$ with the block $D$ projectible
along $V_1$.

 \begin{dfn}
A family of 4D PH-submanifolds $\Phi_\alpha$ intersecting by a PH
curve $\mathcal{C}$ is a {\em pencil\/} if there exists a 2D
PH-foliation in a neighborhood of $\mathcal{C}$ such that the
projection along it is a PH-map and each $\Phi_\alpha$ is projected
to a PH-curve.
 \end{dfn}
In other words in a neighborhood of the curve $J$ is represented by
the above upper-triangular block form. Then for such a pencil the
tensor $N_J$ is degenerate.

Let us recall basics about degenerations of linear Nijenhuis tensors \cite{K$_3$}. Such a tensor can be considered as a $\C$-antilinear map $N_J:\La^2_\C T\to T$ of vector spaces of $\dim_\C=3$. So if $N_J\ne0$ the following situations are possible:
 \begin{itemize}
 \item[] \hspace{-10pt}NDG: $\dim_\C\op{Im}N_J=3$ (non-degenerate);
 \item[] \hspace{-10pt}DG$_1$: $\dim_\C\op{Im}N_J=2$ (weakly degenerate);
 \item[] \hspace{-10pt}DG$_2$: $\dim_\C\op{Im}N_J=1$, there is a kernel $V\in\op{Gr}_1^\C(T)$, $N_J(V,\cdot)=0$.
 \end{itemize}
Generically a pencil belongs to DG$_1$ case. However DG$_2$ can be
obtained in the two following cases: 1. $A$ is projectible along
$V_{23}$ and $B=0$ or $A$ is constant in the above splitting and $B$
is projectible along $V_1$. Then if the tensor $N_J\ne0$, its kernel
coincides with $V_1$. 2. Almost complex structure $D$ on $V_{23}$ is
integrable. Then if the tensor $N_J\ne0$, its kernel is transversal
to $V_1$.

 \begin{prop}\label{prpr}
Let $\Phi_\a$ be a family of 4-dimensional PH-foliations of 6-dimen\-sional $(M,J)$, intersecting by a common foliation $V$ by PH-curves, such that shifts along $V$ is a symmetry of the family as foliations.

Let cardinality of indices $\a$ be at least 4 and at almost every point $x$ there be 4 leaves $\Phi_{\a_i}$ with $T_x\Phi_{\a_i}/T_xV$ of general position in $T_xM/T_xV$.

Then the family is a pencil: There exists a PH-submersion $\pi:(M^6,J)\to(W^4,\ti J)$ with $V$-fibers.
 \end{prop}

Before proving this let us discuss the problem how an almost complex structure $J$ is characterized by its PH-submanifolds. This question is non-void even in dimension 4, on which we concentrate.
In this case the problem can be reformulated as a PH-analog of plane webs.

 \begin{lem}\label{L1}
Let $\Psi_a$ be a PH 4-web of almost complex $(W^4,J)$, i.e. there are foliations $\Psi_a$, $1\le a\le4$, by PH-curves, none two of them being tangent anywhere. Then $J$ is determined by $\Psi_a$ up to sign.
 \end{lem}

 \begin{proof}
We will prove a more general statement: Let $\Psi_a$ be a 4-webs of surfaces in $W^4$ with the same condition of general position at each point. Then there are at most two almost complex structures $\pm J$ making $\Psi_a$ into PH-web.

Indeed, this is the question of linear algebra. We have $T_xW=\Pi_1\oplus\Pi_2$, where $\Pi_a=T_x\Psi_a$ are 2-dimensional subspaces. Complex structures on $\Pi_1$ and $\Pi_2$ determine that on $T_xW$. Since $\Pi_3$ is a complex subspace it is a graph of a complex linear map $F:\Pi_1\to\Pi_2$. This map is nondegenerate and the complex structure on $\Pi_1$ determines that of $\Pi_2$. Now using $\Pi_4$, which is also a graph, we get a complex automorphism $L:\Pi_1\to\Pi_1$, not proportional to identity. So no two different (up to sign) complex structures can commute with it. This proves uniqueness.

Let us discuss existence. It is equivalent to the claim the the spectrum of $L$ is purely complex. Necessity is obvious: if $\op{Sp}(L)$ is real simple or $L$ is a Jordan box, no rotation can commute with it. On the other hand, if $\op{Sp}(L)=\{\frac{\l\pm i}\b\}$, then $J=\b L-\l I$ is a complex structure on $\Pi_1$ and this gives the complex structure on $T_xW$.
 \end{proof}

The above problem is equivalent to the following: Given a family of PH-foliations $\Psi_a$ on $(W^{2n},J)$ and a diffeomorphism $f:W\to W$, mapping them to PH-foliations, how large should be the index set $\{a\}$ to ensure that $f$ is a PH-map or anti-PH: $f^*J=\pm J$. Imposing general position of leaves, making it into PH-web, the modification of the above proof gives the answer $n+2$.

 \medskip

 \begin{Proof}{Proof of Proposition \ref{prpr}}
Shift along transversal $V$ maps transversal foliations $\Phi_\a/V$
into themselves. Since they are complex PH-lines in $TM/TV$, Lemma
\ref{L1} implies that the complex structure in quotient tangent
spaces is preserved. Thus shifts along $V$ preserve the almost
complex structure $J$ in normal direction. Thus the complex
structure becomes of the upper-triangular block form and the result
follows.\qed
 \end{Proof}

 \medskip

Let us call pencils from this Proposition 4-pencils, because there
are 4 foliations in it (but then it extends to a continuous family).

 \begin{rk}
This proposition has certain generalizations to dimensions $2n>6$,
but then one should make more specifications (dimension of
PH-foliations in the pencil, their number etc), so we do not discuss
it.
 \end{rk}

%%%%%%%%%%%%%%%%%%%%%%%%%%%%%%%%%%%%%%%%%%%%%%%%%%%%%%%%%%%%%%%%%%%%%%%%%%%%
%5%
\section{Criteria of integrability}

 \abz
We present several approaches basing on existence of many
PH-submanifolds. It was shown in \cite{K$_2$} that whenever through
every point $x\in(M^{2n},J)$ and every complex
$[\dim_\C=k]$-dimensional subspace in $T_xM$ passes a PH-submanifold
of dimension $2k$ (or PH 2-jet), then $J$ is integrable ($k$ is
fixed).

But with this we require an infinite number of PH-submanifolds to ensure integrability.
This requirement can be much weakened with the same conclusion. We will specify as above to the case $n=3$.

{\bf 1.} We can use the pencils of Proposition \ref{prpr} to get
another criterion as follows. Consider 5 foliations $V_i$ of $M^6$
by PH-curves (these always exist locally), $1\le i\le 5$, none two
of which are tangent and none three have complex dependent tangents
at almost any point.

 \begin{theorem}
Assume that $(M,J)$ admits 10 PH-foliations $\Phi_{ij}$, such that
$\Phi_{ij}$ contain both $V_i$ and $V_j$ and ar symmetric with
respect to shifts along them. Then the structure $J$ is integrable.
 \end{theorem}

 \begin{proof}
Indeed in this case $(M,J)$ admits 5 pencils of PH-foliations of
dimension 4. Proposition \ref{prpr} applies. In fact in this case
the pencils become just as in Example 2 (without upper-triangular
modification) because for any family $A_a=\{\Phi_{ak}\}_{k\ne a}$
there is a transversal PH-foliation $\Phi_{ij}$, $i,j\ne a$. Thus
the tensor $N_J$ is degenerate.

Two pencils can have weak degeneracy along the same complex 2-plane from $\op{Gr}_2(TM,\C)$, but then the next two show another weak degeneracy, so that there is a kernel. The last pencil gives a weak degeneracy of $N_J$, independent of this kernel, whence the Nijenhuis tensor vanishes.
 \end{proof}

The hypotheses of the theorem can be modified to have four 4-pencils, each having 3 common PH-foliations with the other pencils, leading to the same conclusion. However this provides the same total amount 10 of PH-foliations.

{\bf 2.} We can skip organizing PH-foliations in pencils and get the
same claim, but then the number of foliations should grow.

Let us call family $\Phi_\a$ of PH-foliations of dimension 4 quadratically non-degenerate if at almost every point $x\in M$ the tangents $T_x\Phi_\a\in\op{Gr}_2(T_xM,\C)\simeq \C P^2$ do not belong to any real quadric of codimension 1. Note that any 14 points in $\C P^2$ do belong to a real quadric.

 \begin{theorem}
Let $\Phi_\a$ be a family of PH-foliations of dimension 4 in
$(M^6,J)$, $\a=1,\dots,15$. If it is quadratically non-degenerate,
then $J$ is integrable.
 \end{theorem}

 \begin{proof}
If $N_J\ne0$, then the Grassmanian of 4-planes in $T_xM$, which are
invariants with respect to both $J$ and $N_J$ is a real quadric of
codimension 2 in NDG case or codimension 1 in DG cases of
\cite{K$_2$} (it can be also empty, then its codimension is 4). But
no 15 generic points in $\op{Gr}_2^\C(T_xM)$ can belong to a
quadric.
 \end{proof}

{\bf 3.} We can have some intermediate criteria between approach 1,
using fewer number of PH-submanifolds though with some integrability
assumptions, and approach 2, using larger number of PH-submanifolds
but only genericity conditions. For example, assume we have 14
families of 4D PH-foliations of $(M^6,J)$, which have generic
arrangements of tangents at almost every point.

Then we have a field of quadrics $Q_x\subset T_xM$, $x\in M$. If the
structure is non-integrable, this field satisfies certain
integrability criteria ($\Xi(\pi_*\Theta_H(\Pi))=0$ from
\cite{K$_2$}). This is a binding requirement.

 \begin{theorem}
Let $\Phi_\a$ be a family of PH-foliations of dimension 4 in
$(M^6,J)$, $\a=1,\dots,14$ with generic arrangements of tangents
a.e. If the corresponding family of quadrics $Q$ is non-integrable,
then $J$ is integrable.
 \end{theorem}

\small
%%%%%%%%%%%%%%%%%%%%%%%%%%%%%%%%%%%%%%%%%%%%%%%%%%%%%%%%%%%%%%%%%%%%%%%%%%%%

 \vspace{-10pt} \hspace{-20pt} {\hbox to 12cm{ \hrulefill }}
\vspace{-1pt}

{\footnotesize \hspace{-10pt} Institute of Mathematics and
Statistics, University of Troms\o, Troms\o\ 90-37, Norway.

\hspace{-10pt} E-mails: \quad kruglikov\verb"@"math.uit.no} \vspace{-1pt}

\end{document}